\documentclass[a4paper,reqno]{amsart}
\usepackage{amssymb,amscd}

\def\dd#1{\frac{\partial}{\partial#1}}

\newcommand{\Lie}{\mathcal{L}}          
\renewcommand{\div}{\mathrm{div}\,}          

\newtheorem{theorem}{Theorem}
\newtheorem{corollary}{Corollary}
\newtheorem{proposition}{Proposition}
\newtheorem{definition}{Definition}
\newtheorem{lemma}{Lemma}
\newtheorem{example}{Example}
\newtheorem{remark}{Remark}
\begin{document}

\title[Classification of generic multi-vector fields of top
  degree]{Global classification of generic multi-vector fields of top
  degree}

\author{David Mart\'inez Torres}
\address{Departamento de Matem\'aticas, Universidad Carlos III de
 Madrid, Avda. de la Universidad 30, 28911 Legan\'es, Espa\~na}
\email{dmtorres@math.uc3m.es}

\begin{abstract}
For any closed oriented manifold $M$,  the top degree multi-vector fields
transverse to the zero section of $\wedge^{\mbox{\tiny top}}TM$ are classified, up to 
orientation preserving diffeomorphism, in terms of the topology of the arrangement of 
its zero locus and a finite number of numerical invariants. The group governing the 
infinitesimal deformations of such multi-vector fields is computed, and an explicit set
of generators exhibited. For the sphere $S^n$, a correspondence between certain isotopy 
classes of multi-vector fields and classes of weighted signed trees is established.
 \end{abstract}
\maketitle

\section{Introduction}

The recent classification by O.~Radko \cite{Ol} of generic Poisson
structures on oriented surfaces, raises the question of whether it is
possible to extend it to higher dimensions. This classification,
although stated in the language of Poisson geometry, relies on general
results from differential geometry and the classification of area
forms on closed surfaces. The reason is that, in dimension 2, the
integrability condition that a bi-vector field must satisfy in order
to be Poisson is void. So, for generic Poisson structures on an
oriented surface $\Sigma$, the difficult problem of classifying
solutions of a non-linear PDE reduces to the classification of
(generic) sections of a trivial line bundle
$\mathfrak{X}^2(\Sigma)\equiv\Gamma(\wedge^2(T\Sigma))$. Then,
standard methods from differential geometry apply and the problem is
greatly simplified. In this note we show how Radko's classification
can be extended to higher dimensions for generic multi-vector fields
of top degree.

On a surface a bi-vector field defines a Poisson structure. More generally, a
multi-vector field of top degree defines a \emph{Nambu structure} of top
degree. Nambu structures are natural generalizations of
Poisson structures: a \emph{Nambu structure} of degree $r$, on a
manifold $M$, is a $r$-multilinear, skew-symmetric bracket,
\[\{\cdot,\dots,\cdot\}\colon \underbrace{C^\infty
  (M)\times\cdots\times C^\infty (M)}_{r} \rightarrow C^\infty (M),\]
which satisfies the Leibniz rule in each entry, and a Fundamental
Identity \cite{Ta} that naturally extends the Jacobi identity. For top degree
structures the Fundamental Identity is void \cite{Va}.
In spite of the formal similarities between Nambu structures and
Poisson structures, for $r>2$ the Fundamental Identity imposes much
more restrictive conditions than one would expect from the Jacobi
identity. That is, Nambu structures are in a sense harder to find than
Poisson structures. On the other hand, Nambu structures are easier to
describe.

We are interested in generic Nambu structures of top degree on a
closed oriented manifold $M$. By  Leibniz' rule, such a structure is
described by a multi-vector field $\Lambda\in\mathfrak{X}^{\mbox{\scriptsize top}}(M)$.
Genericity means that the graph of $\Lambda$ cuts the zero section of the line
bundle $\wedge^{\mbox{\scriptsize top}} T M$ transversally. In particular, the zero
locus $\mathcal{H}$ of the multi-vector field $\Lambda$ is a
hypersurface in $M$. We will show how one can attach to each connected
component $H^i$ of $\mathcal{H}$ a numerical invariant, called the
\emph{modular period}, which depends only on the germ of $\Lambda$ at
$H^i$. We construct also a global invariant  which measures the ratio
between the volumes of the connected components of the complement of
$\mathcal{H}$, called the \emph{regularized Liouville volume}. These
notions generalize corresponding notions for 2-dimensional Poisson
manifolds.

Our main result is the following:

\begin{theorem}
\label{thm:main}
A generic Nambu structure $\Lambda\in\mathfrak{X}^{\mbox{\scriptsize top}}(M)$ is
determined, up to orientation preserving diffeomorphism, by the diffeomorphism type of the oriented pair $(M,\mathcal{H})$ together with
its modular periods and regularized Liouville volume.
\end{theorem}

For dimension 2 this result recovers the classification of
\cite{Ol}. 

Using Theorem \ref{thm:main} we are also able to describe the Nambu
cohomology group $H^2_{\Lambda}(M)$ which determines the infinitesimal
deformations of the Nambu structure. On the other hand, we show that
for dimension larger than $2$, the Nambu cohomology group
$H^1_{\Lambda}(M)$, which determines the outer automorphisms of the
structure, is infinite dimensional.

The plan of the paper is as follows. In Section 1, we recall the
definition of a Nambu structure of degree $r$ and list briefly some of
its main properties. In Section 2, we consider generic Nambu
structures of degree $n$ on an $n$-dimensional oriented manifold. We define, for
each hypersurface $H$ where the $n$-vector field $\Lambda$ vanishes, a
couple of equivalent invariants. They are the modular $(n-1)$-vector
field $X_\Lambda^H$ and the modular $(n-1)$-form $\Omega_\Lambda^H$,
which give two equivalent ways of describing the linearization of
$\Lambda$ along $H$. In Section 3 we introduce the modular period
$T_\Lambda^H$, which is just the integral (or cohomology class) of the
modular $(n-1)$-form, and depends only on the values of $\Lambda$ on a
tubular neighborhood of $H$. Conversely, we can recover the
Nambu structure on a tubular neighborhood of the oriented
hypersurface $(H,\Omega_\Lambda^H)$ once the modular period
$T_\Lambda^H$ is specified. The proof of theorem \ref{thm:main} is completed  in
Section 4, where we also introduce the regularized Liouville
volume. In Section 5, among the possible cohomological structures one can attach
to a Nambu structure, we consider (i) the group of infinitesimal outer
automorphisms and (ii) the group of infinitesimal deformations of the
structure. The later will turn out to have as many generators as the
numerical invariants above and we will exhibit explicitly a set of
generators, which extends that of 2-dimensional Poisson manifolds. On
the other hand, we will show that the first cohomology group is infinite
dimensional for $n\geq 3$, something to be expected from the local
computation of this group presented in \cite{Mo}. Finally, in
Section 6, we observe that the correspondence between isotopy classes
of generic bi-vectors on
$\Sigma=S^2$ and isomorphism classes of  weighted signed trees given in 
\cite{Ol}, holds for those generic Nambu structures in $S^n$ for
which the zero locus $\mathcal{H}$ only contains spheres.

\subsection*{Acknowledgements}
 I would like to express my deep gratitude to  R. L. Fernandes for
the many valuable discussions that helped to clarify several aspects of the
present work. I would also like to thank A. Serra for his remarks
about the Casimir group, A. Ibort for his many useful comments and the referee for 
her/his valuable suggestions. Finally, I 
would like to thank the Instituto Superior T\'ecnico de Lisboa for
its hospitality during my stay there.

\section{Nambu structures}

Poisson manifolds $(M,\{\cdot, \cdot \})$ are the phase spaces relevant for
Hamiltonian mechanics. For a Hamiltonian system the evolution of any
observable $f \in C^\infty(M)$ is obtained by solving the o.d.e.
\[\frac{df}{dt}=\{H,f\},\]
where $H\in C^\infty(M)$ is the Hamiltonian, a conserved quantity for
the system (the ``energy''). In $1973$ Nambu \cite{Na} proposed a generalization of Hamiltonian
mechanics based on an $n$-ary bracket. The dynamics of an observable $f
\in C^\infty(M)$ would be governed by the analogous o.d.e.
\[\frac{df}{dt}=\{H_1,...,H_{n-1},f\},\]
associated to $n-1$ Hamiltonians $H_1,...,H_{n-1}$, so now we would have
$n-1$ conserved quantities.

In order to have the ``expected'' dynamical properties this bracket
had to satisfy certain constraints. These were clarified by Takhtajan
\cite{Ta}, who gave the following axiomatic definition of a Nambu
structure.

\begin{definition} 
A \textbf{Nambu structure} of degree $r$ in a manifold $M^n$, where
$r\leq n$, is an $r$-multilinear, skew-symmetric bracket,
\[\{\cdot,\dots,\cdot\}\colon \underbrace{C^\infty
  (M)\times\cdots\times C^\infty (M)}_{r} \rightarrow C^\infty (M),\]
satisfying: 
\begin{enumerate} 
\item[(i)] the Leibniz rule:
\[ \{fg,f_1,\dots,f_{r-1}\}=f\{g,f_1,\dots,f_{r-1}\}+\{f,f_1,\dots,f_{r-1}\}g,\]
\item[(ii)] the Fundamental Identity: 
\[\{f_1,\dots,f_{r-1},\{g_1,\dots,g_{r}\}\}=
\sum_{i=1}^{r}\{g_1,\dots,\{f_1,\dots.,f_{r-1},g_i\},\dots,g_r\};\]
\end{enumerate}
\end{definition}

The Liebniz rule shows that the operator
$X_{f_1,\dots,f_{r-1}}:C^\infty (M)\to C^\infty (M)$  which is associated
to $r-1$ functions $f_1,\dots,f_{r-1}$ by
\[ X_{f_1,\dots,f_{r-1}}(g)=\{g,f_1,\dots,f_{r-1}\},\]
is a derivation and hence a vector field. This is called the
\emph{Hamiltonian vector field} associated with $f_1,\dots,f_{r-1}$. More
generally, the Leibniz identity shows that we have an $r$-vector field
$\Lambda\in\mathfrak{X}^r(M)$ such that
\[ \Lambda(df_1\wedge\cdots\wedge df_r)=\{f_1,\dots,f_{r}\}.\]

On the other hand, the Fundamental Identity is
equivalent to the fact that the flow of any Hamiltonian vector field
$X_{f_1,\dots,f_{r-1}}$ is a \emph{canonical transformation}, i.e, 
preserves Nambu brackets. Its infinitesimal version reads
\[ \Lie_{X_{f_1,\dots,f_{r-1}}}\Lambda=0.\]

For  Nambu structures of top degree the Fundamental Identity becomes void \cite{Va}. 

\begin{example}
On $\mathbb{R}^n$ we have a canonical, top degree, Nambu structure
which generalizes the canonical Poisson structure on
$\mathbb{R}^2$. The Nambu bracket assigns to $n$ functions
$f_1,\dots,f_n$ the Jacobian of the map $\mathbb{R}^n\to\mathbb{R}^n$,
$x\mapsto (f_1(x),\dots,f_n(x))$, so that
\[ \{f_1,\dots,f_n\}=\det\left[\frac{\partial f_i}{\partial
    x_j}\right].\]
    
More generally, any volume form $\mu\in\Omega^{\mbox{\scriptsize top}}(M)$ on a
manifold $M$ determines a Nambu structure: if $(x^1,\dots,x^n)$ are
coordinates on $M$, so that $\mu=fdx^1\wedge\cdots\wedge dx^n$, then
the Nambu tensor field is 
\[\Lambda\equiv\frac{1}{\mu}=\frac{1}{f}\frac{\partial}{\partial
  x^1}\wedge\cdots\wedge \frac{\partial}{\partial x^n}.\]
\end{example}

The Fundamental Identity for $r>2$ is of a more restrictive nature
than  the case $r=2$, where it reduces to the
usual Jacobi identity: if $r>2$ besides requiring the fulfilment of a
system of first order quadratic partial differential equations, the
coefficients \emph{must also} satisfy certain system of quadratic
algebraic equations. For example, for a locally constant $r$-vector field 
the system involving first derivatives is
automatically satisfied, while the algebraic relations are
non-trivial, and in fact coincide with the well-known Pl\"ucker
equations. Hence, only decomposable $r$-vectors define constant Nambu
structures. Another example of this rigidity is the following well-known
 proposition (see \cite{Ta}):

\begin{proposition}
Let $\Lambda$ be a Nambu structure. For any function
$f\in C^\infty(M)$, the contraction $i_{df}\Lambda$ is also a Nambu
structure. 
\end{proposition}

This rigidity makes it harder to ``find'' Nambu structures than
Poisson structures. On the other hand, it makes Nambu structures
easier to describe. Henceforth, we will assume that $r>2$ if $n\ge 3$.

First of all, the Hamiltonian vector fields span a generalized
foliation for which the leaves are either points, called
\emph{singular points}, or have dimension equal to the degree of the
structure. Around these \emph{regular points} we have the following
canonical form for a Nambu structure (see for example \cite{Va}):

\begin{proposition}
Let $x_0\in M$ be a regular point of a Nambu structure $\Lambda$ of
degree $r$. There exist local coordinates $(x^1,\dots,x^n)$ centered at
$x_0$, such that
\[ \Lambda=\frac{\partial}{\partial x^1}\wedge\cdots\wedge
\frac{\partial}{\partial x^r}.\] 
\end{proposition}

For the singular points there are some deep linearization results due
to Dufour and Zung \cite{DZ}.

\section{Generic Nambu structures of top degree}
In this section we consider Nambu structures of degree $n$ in a
closed orientable $n$-dimensional manifold $M$. Notice that in this
case the Fundamental Identity is void \cite{Va}, so a Nambu structure is just a
multi-vector field $\Lambda\in\mathfrak{X}^n(M)$. We will restrict our
attention to \emph{generic} Nambu structures:

\begin{definition} 
A Nambu structure $\Lambda\in\mathfrak{X}^n(M)$ is called
\textbf{generic} if its graph cuts the zero section of the line bundle
$\wedge^n TM$ transversally.
\end{definition}

The generic sections form an open dense set in the Whitney $C^\infty $
topology. 

Let us fix a generic $\Lambda\in\mathfrak{X}^n(M)$.
Its set of zeros, denoted $\mathcal{H}$, is the union of a finite
number of connected hypersurfaces: $\mathcal{H}=\bigcup_{i \in I}H^i$, $\#I<\infty$.
Fix one of them and call it $H$. 

Over the points of $H$ there is some linear information attached to
$\Lambda $, namely the \emph{intrinsic derivative} $d\Lambda^H \in
T^*_HM\otimes \wedge^nTM$. It can be defined as
$d\Lambda^H\equiv\nabla\Lambda_{\mid H}$, where $\nabla$ is any linear
connection on $\wedge^nTM$. This is independent of the choice of
connection. The intrinsic derivative gives the linearization of the
Nambu structure at $H$: if we view $\Lambda$ as a section, it is the
tangent space to the graph of $\Lambda$. It is important to observe that
$d\Lambda^H$ never vanishes due to the transversality assumption.
Notice that $d\Lambda^H$ is a section of $T^*_HM\otimes \wedge^nTM$, but
due to the nature of our (trivial) line bundle it has two equivalent
interpretations which we shall now explain.

Fix a volume form $\Omega$ in some neighborhood of $H$ in
$M$, so that after contracting it with the factor of $d\Lambda^H$ in $\wedge^nTM$, 
we obtain $d\Lambda^H\otimes \Omega \in T^*_HM$. 

\begin{definition} 
The \textbf{modular $(n-1)$-vector field} of $\Lambda$ along $H$ 
is the unique $(n-1)$-vector field
$X_{\Lambda}^H\in\mathfrak{X}^{n-1}(H)$ such that
$i_{X_{\Lambda}^H}\Omega=d\Lambda^H \otimes\Omega$. 
\end{definition}

This definition does not depend on the choice of $\Omega$: if
$\tilde{\Omega}$ is another volume form, then 
$\tilde{\Omega}=f\Omega$ for some non-vanishing smooth function $f$,
and we find 
\[i_{X_\Lambda^H}f \Omega= d\Lambda^H \otimes f\Omega.\]

Notice that since $X_{\Lambda}^H$ is ``tangent'' to $H$ (i.e. $X_{\Lambda}^H\in \mathfrak{X}^{n-1}(H)$) and no-where vanishing,  we can define the \emph{modular $(n-1)$-form} along $H$ to be the dual $(n-1)$-form
$\Omega_{\Lambda}^H\in \Omega^{n-1}(H)$; that is, $\Omega_{\Lambda}^H(X_{\Lambda}^H)=1$.
If we fix a vector field $Y$ over $H$, which is transverse to $H$, the modular 
form along $H$ is given by
\[
\Omega_\Lambda^H=(-1)^{n-1}\frac{1}{i_Yd\Lambda^H\otimes\Omega}j^*i_Y\Omega,
\] 
where $j\colon H\hookrightarrow M$ is the inclusion and 
$i_Yd\Lambda^H\otimes\Omega\in C^\infty(H)$ the no-where vanishing function 
obtained by contracting $Y$ with $d\Lambda^H\otimes\Omega$. This expression
is independent of $Y$. The modular $(n-1)$-form along $H$ is everywhere non-zero,
and hence $\Lambda$ determines an orientation in $H$. It is clear that
to give either one of $d\Lambda^H$, $X_\Lambda^H$ or $\Omega_\Lambda^H$,
determines the others.

Let us relate these definitions with the well-known notion of modular class
of a Poisson manifold. For any Nambu structure of degree $r$ on an
oriented manifold there is a natural generalization of the modular
class of a Poisson manifold \cite{M}, which we now recall. Again we fix a
volume form $\Omega$ on $M$. Then, for any $r-1$ functions
$f_1,\dots,f_{r-1}$ on $M$, we can compute the divergence of the
corresponding Hamiltonian vector field:
\[ 
(f_1,\dots,f_{r-1})\mapsto \div^{\Omega}(X_{f_1,\dots,f_{r-1}})\equiv\frac{1}{\Omega}\Lie_{X_{f_1,\dots,f_{r-1}}}\Omega.\]

It turns out that this assignment defines a $(r-1)$-vector field
$\mathcal{M}_\Lambda^{\Omega}$ on $M$. If $\tilde{\Omega}=g\Omega$ is
another volume form, where $g$ is some no-where vanishing smooth function,
we have
\[\mathcal{M}_\Lambda^{\tilde{\Omega}}=\mathcal{M}_\Lambda^{\Omega}+X_g,\]
where $X_g$ is the $(r-1)$-vector field
\[ X_g(f_1,\dots,f_{r-1})=\{f_1,\dots,f_{r-1},g\}.\]

One can introduce certain Nambu cohomology groups to take care of this
ambiguity so that the cohomology class
$[\mathcal{M}_\Lambda^{\Omega}]$ is well-defined and independent of
$\Omega$. This class is called the \emph{modular class} of the Nambu
manifold $M$ and is the obstruction for the existence of a volume form
on $M$ invariant under Hamiltonian automorphisms.

Back to top degree Nambu structures, each volume form $\Omega$ determines
a modular $(n-1)$-vector field 
$\mathcal{M}_\Lambda^{\Omega}$ representing the modular class, and
which will depend on $\Omega$. However, at points where the Nambu
tensor vanishes all modular vector fields give the same value (see
\cite{M}). The modular $(n-1)$-vector field $X_\Lambda^H$ along $H$,
that we have introduced above, is nothing but the restriction of any
modular vector field to $H$. In our case, however, it has the
additional properties that it is non-zero and tangent to $H$.

\section{Local characterization of a Nambu structure}

In this section we study the local behavior of a generic Nambu structure
$\Lambda\in\mathfrak{X}^{n}(M)$ on a neighborhood of its zero
locus. We show that the germ of $\Lambda$ around a connected
component $H$ of its zero locus is determined, up to isotopy, by the
modular periods (to be introduced below).

Since the intrinsic derivative is functorial we immediately conclude
that

\begin{lemma}
Given two generic Nambu structures $\Lambda_1$ and $\Lambda_2$ with
zero locus $\mathcal{H}_1$ and $\mathcal{H}_2$, and a diffeomorphism of
Nambu structures $\psi\colon
(M,\Lambda_1)\longrightarrow(M,\Lambda_2)$ then
\[ \psi_*X_{\Lambda_1}^{\mathcal{H}_1}= X_{\Lambda_2}^{\mathcal{H}_2},
{\mbox{and}} \ \psi_*\Omega_{\Lambda_1}^{\mathcal{H}_1}=
\Omega_{\Lambda_2}^{\mathcal{H}_2}.\]
\end{lemma} 

Hence, it follows that a necessary condition for such a map to exist
 is that the cohomology classes
$[\Omega_{\Lambda_2}^{\mathcal{H}_2}]$ and
$[\Omega_{\Lambda_1}^{\mathcal{H}_1}]$ correspond to each other.

Now recall that, given a generic Nambu structure $\Lambda$, each
connected component $H$ of its zero locus $\mathcal{H}$ has an induced
orientation from the $n$-vector field $\Lambda$. Hence, a class in
$H_{dR}^{n-1}(H)$ is completely determined by its value on the
fundamental cycle $H$.

\begin{definition} 
The \textbf{modular period} $T_\Lambda^H$ of the component $H$ of the
zero locus of $\Lambda$ is 
\[ T_\Lambda^H\equiv\int_H{\Omega_\Lambda^H}>0.\]
\end{definition}

In fact, this positive number determines the Nambu structure in a
neighborhood of $H$ up to isotopy. To prove that we need the
following classical result concerning the classification of volume forms.

\begin{lemma}[(Moser)] 
Let $M$ be an orientable closed manifold, $\Omega_1$ and $\Omega_2$ two
volume forms in $M$. If
$[\Omega_1]=[\Omega_2]\in H^{\mbox{\scriptsize{top}}}(M)$, there exists a
diffeomorphism isotopic to the identity which sends $\Omega_1$ to $\Omega_2$.
Moreover, it can be chosen to have support in the closure of the complement 
of the closed set where the two volume forms coincide.
\end{lemma}

The above result can be adapted to volume forms in compact manifolds with
boundary which coincide in neighborhoods of the boundary components.
 We can now state and prove the main result in this section.

\begin{proposition}\label{prop:local}
Let $\Lambda_1$ and $\Lambda_2$ be generic Nambu structures in $M$ (oriented) which
share a common component $H$ of their zero locus with equal induced orientation, 
and for which the modular periods coincide: $T_{\Lambda_1}^H=T_{\Lambda_2}^H$.
Then, there exists a diffeomorphism $\varphi:M\to M$, isotopic to the identity,
and neighborhoods $U_1$ and $U_2$ of $H$, such that $\varphi$ sends
$(U_1,\Lambda_1)$ to $(U_2,\Lambda_2)$.
\end{proposition}

\begin{proof}
First we can use Moser's lemma to construct a diffeomorphism $\phi:
M\to M$ isotopic to the identity, which maps $H$ to itself and sends
$\Omega^H_{\Lambda_1}$ to $\Omega^H_{\Lambda_2}$. Hence, we can
assume that $\Omega^H_{\Lambda_1}=\Omega^H_{\Lambda_2}$, and the
problem reduces to a global linearization one.

We fix a collar $U=[-1,1]\times H$ of the hypersurface $H$, with transverse
coordinate $r$. Denoting the Nambu structure by $\Lambda$ we define
$\Lambda_0=(-1)^{n-1}\dd r \wedge X_{\Lambda}^H $. We can write
$\Lambda=f\Lambda_0 $ for some $f\in C^\infty(U)$ and the linearization of 
$\Lambda$ is $\Lambda_1=r\Lambda_0$. The o.d.e. for a change  of coordinates 
$\phi(r,x)=(g(r,x),x)$ ($x$ a coordinate in $H$) so that $\phi_*f\Lambda_0=r\Lambda_0$ is:

\begin{equation}\label{eq:chcoord}
\frac{\partial g}{\partial r}=\frac{g}{f}\,
\end{equation}

 Since $f$ has for each $x$ an expansion along the radial direction of the form $r+a_2(x)r^2+\cdots $, 
 any solution  $g(r,x)=ke^{\int{\frac{1}{f}dr}}$, $k>0$, defines a smooth 
orientation preserving change of coordinates fixing $H$.
\end{proof}

\begin{remark}
The existence of a one parameter family of solutions for equation \ref{eq:chcoord}
reflects the fact that for any linear Nambu structure $cr\dd r \wedge
X_{\Lambda}^H$, $c\in \mathbb{R}\smallsetminus\{0\}$, rescaling the radial coordinate 
is a transformation preserving the Nambu structure (canonical transformation). 
Note also that the reflection along $H$ is a canonical transformation, but it changes 
the orientation of the tubular neighborhood.  We are interested in
transformations isotopic to the identity so our choice of $g$ above is with $k>0$.
\end{remark}

\section{The global description of Nambu structures}

In order to have a diffeomorphism between two Nambu manifolds it is
necessary to have a diffeomorphism sending the zero locus of one
structure to the zero locus of the other, preserving their induced
orientations. Assuming this condition to hold, our problem is that of
transforming a generic structure $\Lambda_1$ into another $\Lambda_2$,
with common oriented zero set $\mathcal{H}=\bigcup_{i \in I}H^i$.

First of all, in the previous section we proved that if the modular
periods of each component coincide, we can find collars $U^i_1$ and
$U_2^i$ of the hypersurfaces $H^i$ and a diffeomorphism isotopic to
the identity $\varphi$ sending  $(\mathcal{U}_1,\Lambda_1)$ to
$(\mathcal{U}_2,\Lambda_2)$, where $\mathcal{U}_j= \bigcup_{i \in I}U_j^i$, $j=1,2$.

 The restriction of the Nambu structures to the connected components of 
 $M\smallsetminus\mathcal{H}$ define volume forms
(the dual volume forms). However, their volumes are infinite so
one cannot require them to match. Instead, we could try to define
finite ratios of the volumes between the various components. 
This raises some accounting problems, so instead we observe that
for a component $H$, a volume form $\Omega$ defined in a neighborhood
of $H$ and the volume form $\Lambda^{-1}$ define orientations on the
complement of $H$, which match on one side of $H$ and are opposite on
the other side. Given any function $h\in C^{\infty}(M)$ only vanishing
linearly in the components of $\mathcal{H}$ (its graph is transverse to the zero
section and vanishes exactly at $\mathcal{H}$), we let
$M^\epsilon(h)=h^{-1}(\mathbb{R}\smallsetminus(-\epsilon,\epsilon))$,
with $\epsilon>0$ small enough so that $M^\epsilon(h)$ contains
the complement of the union of collars of the $H^i$, and we set
\[ V_\Lambda ^\epsilon (h)=\int_{M^\epsilon (h)}\Lambda^{-1}.\] 

Here $\Lambda^{-1}$ denotes the volume form dual to $\Lambda$, and
 to integrate we use the given orientation of $M$. The following
definition generalizes the one given in \cite{Ol} for the case of
2-dimensional Poisson manifolds.

\begin{definition}
\label{reglvol} 
The \textbf{regularized Liouville volume} of $\Lambda $ is defined as
\[ V_\Lambda=\lim_{\epsilon \rightarrow 0}V^\epsilon_\Lambda(h),\]
where $h$ is any function only vanishing linearly at $\mathcal{H}$.
\end{definition}

We only need to
prove the independence on the choice of function $h$, because for a 
function that locally coincides with a radial coordinate in which the 
$n$-vector field is linear the existence of the limit is straightforward. 

Independence on the choice of $h$ can be checked by  adapting Radko's 
proof for surfaces \cite{Ol} to our higher dimensional setting. We fix coordinates
$(r,x)$ around each component $H$ such that $\Lambda=(-1)^{n-1}r\dd r\wedge
X^H_{\Lambda}$, and we consider two cases:
\begin{enumerate}
\item Assume $h(r,x)=g(x)r$, with $g(x)\neq 0$
for all $x\in H$. Let us denote by $H^{>1}(|g|)$ (resp. $H^{<1}(|g|)$) the 
points of $H$ where the absolute value of $g$ is greater (resp. smaller) than $1$. Then
\begin{eqnarray*}
V_\Lambda^\epsilon(h)-V_\Lambda^\epsilon(r)&=&\pm\int_{H^{>1}(|g|)}\left(\int_{[-\epsilon,-\epsilon/|g(x)|]\cup [\epsilon/|g(x)|,\epsilon] }(-1)^{n-1}\frac{dr}{r}\right)\Omega^H_\Lambda \mp\\
&\mp &\int_{H^{<1}(|g|)}\left(\int_{[-\epsilon/|g(x)|,-\epsilon]\cup [\epsilon,\epsilon/|g(x)|]}(-1)^{n-1}\frac{dr}{r}\right)\Omega^H_\Lambda,
\end{eqnarray*}

and each summand vanishes (for every $\epsilon>0$ small enough). 
\item Since $h$ vanishes linearly at $\mathcal{H}$, the function
$h-\frac{\partial h}{\partial r}(0,x)r$ vanishes in the radial
direction at least to second order at $H$. The compactness of $H$  implies
that for all $x\in H$ and all $\epsilon>0$ small enough, positive 
constants $k_1$ and $k_2$ exists (independent of $\epsilon$ and $x$), 
with $a_\epsilon>k_1 \epsilon$ and $b_\epsilon-a_\epsilon<k_2 \epsilon^2$, such that:

\[|V_\Lambda^\epsilon(h)-V_\Lambda^\epsilon({\textstyle\frac{\partial h}{\partial r}(0,x)r)}|\leq \int_H\left(\int_{[-b_\epsilon,-a_\epsilon]\cup [a_\epsilon,b_\epsilon]}\left| \frac{1}{r}\right|dr\right)\Omega_\Lambda^H.\]

Thus
\[|V_\Lambda^\epsilon(h)-V_\Lambda^\epsilon({\textstyle\frac{\partial h}{\partial r}(0,x)r)}|\leq \int_H k\epsilon \Omega_\Lambda^H=k\epsilon T_\Lambda^H,\]

for some $k>0$ and hence when 
$\epsilon\rightarrow 0$ the difference vanishes.
\end{enumerate}

The modular periods and the regularized volume determine the Nambu
structure:
\begin{theorem}\label{main} 
For $j=1,2$, let $M_j$ be oriented closed manifolds with generic
Nambu structures $\Lambda_j$ having zero locus $\mathcal{H}_j=\bigcup_{i \in
I}H_j^i$. Assume that there exists a diffeomorphism $\psi$ sending
$(M_1,\mathcal{H}_1)$ to $(M_2,\mathcal{H}_2)$ and preserving the induced
orientations of the zero locus. Then there exists an isomorphism
between the two  Nambu structures, isotopic to $\psi$, if and
only if the following conditions are satisfied:
\begin{enumerate}
\item[(i)] the modular periods coincide, i.e,
  $T_{\Lambda_1}^{H_1^i}=T_{\Lambda_2}^{\psi H_1^i}, \; \forall i \in I$,
\item[(ii)] the regularized volumes match, i.e,
  $V_{\Lambda_1}=\varepsilon V_{\Lambda_2}$, where $\varepsilon=1$
  if $\psi $ is orientation preserving and $\varepsilon=-1$ if it
  reverses the orientations on the $M_i$. 
\end{enumerate}
\end{theorem}

\begin{proof}
Because of proposition \ref{prop:local} if the modular periods of each component coincide, we
can find collars $U^i_1$ and $U^i_2$ of the hypersurfaces $H^i$, and a
diffeomorphism isotopic to the identity $\varphi$ which sends
$(\mathcal{U}_1,\Lambda_1)$ to $(\mathcal{U}_2,\Lambda_2)$, where $\mathcal{U}_j=
\bigcup_{i \in I}U_j^i$.

The volume of each connected component of $M\smallsetminus\mathcal{H}$ with respect 
to the duals $\Lambda_i^{-1}$ of any of the
$n$-vectors is infinite. For any such connected component $L$, we can select a
hypersurface $H^{i_0}$ on its boundary and shrink accordingly the size
of either $U_1^{i_0}$ or $U_2^{i_0}$ (recall we have canonical
transformations doing that) such that one can find compact submanifolds
$W_j\subset L$  which are the result
of removing from $L$ the corresponding side of the collars of radius $1/2$ say
(the original radius is $1$), verifying:
\begin{itemize}
\item[(i)] the $\Lambda_1^{-1}$-volume of $W_1$ coincides with the $\Lambda_2^{-1}$-volume of $W_2$.
\item[(ii)] $\varphi$ sends $W_1$ to $W_2$. 
\end{itemize}

Finally, we apply Moser theorem
to conclude the existence of a
diffeomorphism isotopic to the identity also matching the Nambu
structures at $L$.

Observe that when we modify the size of $U^{i_0}_1$ say, we are
changing the volume of both $L\smallsetminus\mathcal{U}_1 $ and $L'\smallsetminus\mathcal{U}_1$,
where $L$ and $L'$ are the connected components of $M\smallsetminus\mathcal{H}$
whose boundary contains $H^{i_0}$. It follows that without further assumptions we can make our $n$-vector 
fields coincide, as well as in the collars of $\mathcal{H}$, in all the 
connected components but possibly one. To see this we can take the graph dual
to the splitting given by $\mathcal{H}$, where each vertex represents a connected
component of $M\smallsetminus\mathcal{H}$ and an edge joining two vertices stands
for a connected hypersurface on its common boundary, and consider a maximal tree
(a contractible subgraph containing all the vertices, which always exists); 
we fix a vertex $v_0$ on this tree and consider the graph distance (of vertices)
with respect to $v_0$.  We can then proceed by
stages, where at each stage we consider all the vertices at the same
distance of $v_0$, starting from the furthest way vertices. For those vertices,
i.e, connected components of $M\smallsetminus\mathcal{H}$, we apply the above 
reasoning to the hypersurface representing the only edge reaching them 
(if a furthest way vertex is reached by two different edges, we could join 
the two vertices in the boundary of the edges with $v_0$ and then form a loop).
When we are done we erase those vertices and edges connecting to them, so we
obtain a smaller tree. We keep on doing that until we reach the vertices at 
distance one (notice that each step  does not affect the connected components
corresponding to further vertices, since that would imply the existence of a 
loop in the tree).  The same process is applied to all but one edge. The
fact that the regularized volumes match, grants us the matching of the
volumes of the two remaining components of $M\smallsetminus \mathcal{H}$ for both
Nambu structures, once an  appropriate collar of the hypersurface representing
the last edge has been removed, and this finishes the proof.
\end{proof}

The set of generic Nambu structures has an action of ${\mbox {Diff}_0}(M)$ 
(resp. ${\mbox {Diff}^+}(M)$).  Its space of orbits has as many connected 
components as isotopy classes (resp. oriented diffeomorphism classes) of 
oriented hypersurfaces  $\mathcal{H}=\bigcup_{i\in I}H^i$. Theorem \ref{main} 
gives an explicit parametrization of each connected component of this moduli space.

\section{Nambu cohomology}

There are several cohomology theories one can associate to a Nambu
manifold (see \cite{M,Mo}). Here we will be interested in the
cohomology associated with the complex
\[0\longrightarrow\wedge^{n-1}C^\infty(M) \longrightarrow \mathfrak{X}(M)
\longrightarrow \mathfrak{X}^n(M)\longrightarrow 0 \,,\] 
where the first map is
$f_1\wedge\cdots\wedge f_{n-1}\mapsto X_{f_1,\dots,f_{n-1}}$, while the
second map is $X\mapsto \Lie_X\Lambda$. Notice that the associated
cohomology groups have simple geometrical meanings:
\begin{enumerate}
\item[(i)] $H^0_\Lambda (M)$ is the space of Casimirs of the Nambu structure;
\item[(ii)] $H^1_\Lambda (M)$ is the space of infinitesimal outer automorphisms of
  the Nambu structure;
\item[(iii)] $H^2_\Lambda (M)$ is the space of infinitesimal deformations of
  the Nambu structure.
\end{enumerate}

Computations of Nambu cohomology for \emph{germs} of Nambu structures
defined by quasipolynomials (functions vanishing at the origin with
finite codimension) were done by Monnier in \cite{Mo}. Here we are
interested in \emph{global} Nambu structures with the simplest
singularity. These computations can be thought of as a infinitesimal
version on the classification theorem; in particular, $H^2_\Lambda (M)$ 
will turn out to be the tangent space of the class of $[\Lambda]$ in the 
moduli space of generic Nambu structures.

The main result of this section is the following

\begin{theorem}
Let $\Lambda$ be a generic Nambu structure on an oriented closed
manifold $M$ with zero locus $\mathcal{H}=\bigcup_{i \in I}H^i$.
The group $H^2_\Lambda(M)$ has dimension $\#I+1$ and a set of
generators is given by
\[ \beta_1(-1)^{n-1}r\dd r \wedge X_\Lambda^{H^1},\dots,\beta_{\#I}(-1)^{n-1}r\dd r \wedge  X_\Lambda^{H^{\#I}},\ \Omega,\]
where $\Omega $ is a volume form, and each $\beta_i$ is a bump function
supported in a collar of the hypersurface $H^i$.
\end{theorem}

We can give a geometric description of the isomorphism
$H^2_\Lambda(M)\simeq\mathbb{R}^{\#I+1}$ as follows. Each $\Theta \in
\mathfrak{X}^n(M)$ is cohomologous to an $n$-vector field whose
vanishing set contains $\mathcal{H}$ and is generic in a neighborhood
of $\mathcal{H}$. Then we can write $[\Theta]=[g\Lambda]$ where $g$ is
some smooth function which assumes a constant value $c_i$ in the
collar of each $U^i$. The isomorphism is
\[
[\Theta] \longmapsto
\left(\frac{T_\Lambda^{H^1}}{T_\Theta^{H^1}},\dots,
\frac{T_\Lambda^{H^{\#I}}}{T_\Theta^{H^{\#I}}},V^{\mathcal{H},\Lambda}_\Theta
\right),\]
where:
\begin{enumerate}
\item[(a)] $T_\Lambda^{H^i}/T_\Theta^{H^i}=c_i$,
\item[(b)] $V^{\mathcal{H},\Lambda}_\Theta$ is the regularized
  integral of $g\frac{1}{\Lambda}$.
\end{enumerate}

The rest of this section is dedicated to the proof of this result,
which consists of a Mayer-Vietoris argument: we first compute the
groups in the collars and then we glue them using information about
the infinitesimal automorphisms in those neighborhoods.

\subsection{Computation of $H_\Lambda^2(U)$}
Let us fix $H \subset\mathcal{H}$ and $U=(-1,1)\times H$ a collar in 
which $\Lambda$ is linear.

\begin{proposition} 
$H_\Lambda^2(U)\simeq \mathbb{R}$ and a generator is given by the
linearization $(-1)^{n-1}r \dd r\wedge X_\Lambda^H $.
\end{proposition}

\begin{proof}
Any vector field $X$ can be written $X=A\dd r +X_H$, where $A \in
C^\infty(U)$,  $X_H \in  \pi^*_2TH$, $\pi_2\colon (-1,1)\times H\rightarrow H$ the
projection onto the second factor. Defining $\Lambda_0=(-1)^{n-1}\dd r
\wedge X_\Lambda^H$, one has:
\begin{eqnarray}\label{eq:Lieder}
\Lie_X\Lambda &= & \Lie_X r\Lambda_0=A\Lambda_0+r\Lie_X\Lambda_0=\nonumber\\
           &= & \left(A-r\frac{\partial A}{\partial r}\right)\Lambda_0+(-1)^{n-1}r\dd
r\wedge \Lie_{X_H}X_\Lambda ^H=\nonumber\\
&=&\left(A-r\frac{\partial A}{\partial
r}+r\div^{\Omega_\Lambda^H}(X_H)\right)\Lambda_0,
 \end{eqnarray}
where $\div^{\Omega_\Lambda^H}(X_H)$ is the divergence of $X_H$ with respect to
$\Omega_\Lambda^H $ (for each $r$).

Any $n$-vector field $f\Lambda_0$ in $U$ is
equivalent to a unique linear one: we observe that if $f$ is at least quadratically 
vanishing at the origin, it is a coboundary;  write $f=r^2g$ and  apply equation 
\ref{eq:Lieder} to obtain 

\[\Lie_{(-r\int{gdr})\frac{\partial}{\partial r}}\Lambda=
(-r\int{gdr}+r\int{gdr}+ r^2g) \Lambda_0=f\Lambda_0.\]

We first make $f$ vanish 
at the origin by adding the coboundary $\Lie_{-f\dd r}\Lambda$ (so it becomes
$r\frac{\partial f}{\partial r}\Lambda_0$).  Writing $f=r\hat{f}$, the coefficient
of the linear representative we look for is $c=\int_{\{0\}\times
H}{\hat{f}}\Omega_\Lambda^H \in \mathbb{R}$, because if  we choose any $Y\in\mathfrak{X}(H)$
such that $\hat{f}_{\mid H}+div^{\Omega_\Lambda^H}(Y)=c$, from equation \ref{eq:Lieder}
again we deduce that  $\Lambda_1=f\Lambda_0+\Lie_Y\Lambda$ has $cr\Lambda_0$ as constant linear
part at $H$. Finally, $\Lambda_1-cr\Lambda_0$ is at least quadratically
vanishing at $H$ and hence it is a coboundary.

The uniqueness of the linear representative will follow from the non existence of 
solutions for the equation $E_c\equiv A-r\frac{\partial A}{\partial r}-r\div^{\Omega_\Lambda^H}(X_H)=cr$,
for $c=1$. We need also to study the equation $E_0$
of Nambu infinitesimal automorphisms.

\begin{lemma} 
\label{lemma:infinit}
The equation $E_1$ has no solutions, and $Z^1_\Lambda (U)$ --the
space  of solutions of $E_0$-- can be identified with the vector space:
\[Z^1_\Lambda (U)\cong span\langle r\dd r, X_H\in  \pi^*_2TH \,|\,
\div^{\Omega_\Lambda^H}(X_H(0))=0\rangle.\]
\end{lemma}

\begin{proof}[of Lemma \ref{lemma:infinit}] 
In the equations $E_c$ we can write the term
$\div^{\Omega_\Lambda^H}(X_H)$ in the form $\psi_r$, where $\psi_r$ is a smooth
family in $r\in (-1,1)$ of functions in $C^\infty(H)$ which satisfy
$\int_H{\psi_r \Omega_\Lambda^H}=0,\;\forall r \in (-1,1)$.  The solutions of 
$E_1$ can be explicitly written as: 

\[ A=kr+r\int{\frac{\psi_r-1}{r}dr},\; k\in \mathbb{R}\]

Any solution has to be a smooth continuation of the above  expression, but it
cannot exist. Indeed, since $\psi_0$ has vanishing integral, we can find a point
$x$ in $H$, such that $\psi_0(x)=0$. Hence, in a small segment $[-\epsilon,
\epsilon]\times \{x\}$ the real valued function
$r\int{\frac{\psi_r(x)-1}{r}dr}$ is, up to a smooth function, $r\log r$ (not even
$C^1$).
\end{proof}

With this we have finished the computation of $H^2_\Lambda(U)$.  
\end{proof}

\begin{remark}
Regarding $E_0$, its solutions are of the form
\[A=kr+r\int{\frac{\psi_r}{r}dr},\]
which will be smooth if and only  if $\psi_0=0$, or the corresponding vector
field $X_H(0)$ is divergence free with respect to $\Omega_\Lambda^H$.
Hence,
\[Z^1_\Lambda(U) \cong span<r\dd r, X_H \in  \pi^*_2TH\,|
\,\div^{\Omega_\Lambda^H}(X_H(0))=0 >\]
\end{remark}

\subsection{From $H_\Lambda^2(U)$ to $H_\Lambda^2(M)$}
The remaining step is to piece all the local information. We just showed that in
the same radial coordinate in which $\Lambda $ is linearized in
$U^i=(-1,1)\times H^i$, we can find a representant $\Theta $ of the cohomology
class such that $\Theta_{\mid U^i}=c_i(-1)^{n-1}r\dd r \wedge
X_\Lambda^{H^i}=c_i\Lambda $, where the relative period
$\frac{T_\Lambda^{H^i}}{T_\Theta^{H^i}}$ is $c_i$, which might
be zero. In particular, for a no-where vanishing Nambu structure all the local invariants
vanish because by looking at its dual form it is clear that since it does not
vanishes, we can push its graph down (or up) to the zero section to make it
vanish in the $U^i$'s. This operation can be made without changing the area (do it
randomly and multiply by the ratio of both areas). We also saw that we can
restrict our attention to coboundaries $X$ such that $X_{\mid U^i}$ is a
solution of $E_0$ in the
radial coordinates. 

The global regularized volume with respect to $\Lambda $ is well defined because 
the regularized volume of $\Lie_X\Lambda $ vanishes, for $X$ infinitesimal automorphism
of $\Lambda $ in $U^i$. To see that we choose  a function $h$ coinciding with the 
radial coordinate in each $U^i$, $M^r(h)=M\smallsetminus\bigcup_{i \in I}(-r,r)\times H^i$. We have:

 \begin{equation}\int_{M^r(h)}{\Lie_X\Lambda}=\int_{M^r(h)}{di_X\frac{1}{\Lambda}}=  \pm
\sum_{i\in I}\left(\int_{\{ r \}\times H^i}
{i_X\frac{1}{\Lambda}}-\int_{\{-r \}\times H^i} {i_X\frac{1}{\Lambda}}\right)
\end{equation}

And for a fixed component $H$  and
$X=kr+r\int{\frac{\div^{\Omega_\Lambda^H}(X_H)}{r}dr}\dd r +X_H$, the function

$I(r)=\int_{\{r\}\times H}{i_X\frac{1}{\Lambda}}$ equals:

\begin{eqnarray*}
I(r)&= &(-1)^{n-1}\int_{\{r\}\times H}{\frac{1}{r}
\left(kr+r\int{\frac{\div^{\Omega_\Lambda^H}(X_H)}{r}dr}\right)\Omega_\Lambda^H }=\\
&=&(-1)^{n-1} \int_{\{r\}\times H}{\left(k+\int{\frac{\div^{\Omega_\Lambda^H}(X_H)}{r}dr}\right)\Omega_\Lambda^H }
\end{eqnarray*}

Due to the fact that $\div^{\Omega_\Lambda^H}(X_H)(0)=0$, the above formula  defines a smooth
function  for all $r \in [-1,1]$.  Its derivative is easily computed:
\begin{eqnarray*}
\frac{dI}{dr}&=&(-1)^{n-1}\frac{d}{dr}\int_{\{r\}\times H}{\left(k+\int{\frac{\div^{\Omega_\Lambda^H}(X_H)}{r}dr}\right)\Omega_\Lambda^H}=\\
&=& (-1)^{n-1}\int_{\{r\}\times H}{\frac{\div^{\Omega_\Lambda^H}(X_H)}{r}\Omega_\Lambda^H=0}
\end{eqnarray*}
The vanishing is clear for $r\neq 0$ and follows by continuity. Hence $I(r)$ is constant and
$V^{\mathcal{H},\Lambda}_\Theta $ is well defined.

It only remains to show that two $n$-vectors $\Theta_1$ and $\Theta_2$ with equal linearizations
and regularized volume are in the same class. Its difference has a representative $\tilde{\Theta}$
vanishing in a neighborhood of the boundary of $M-\stackrel{\circ}{\mathcal{U}}$. Then
the form $\frac{1}{\Lambda}(\tilde{\Theta})\cdot\frac{1}{\Lambda} $ has compact support
(shrinking a bit the collars if necessary) and vanishing integral, so we can find
a compactly supported vector field $Y$ whose divergence is
$\frac{1}{\Lambda}(\tilde{\Theta})\cdot\frac{1}{\Lambda} $. It follows that
$\Lie_{\tilde{Y}} \Lambda =\tilde{\Theta}$, where $\tilde{Y}$ extends $Y$ trivially. 

The assertion about the basis of $H^2_\Lambda(M)$ follows easily.

\subsection{Some comments about $H^1_\Lambda(M)$ and $H^0_\Lambda (M)$}

We will focus our attention in what happens in a collar $U$. Given
$f \in C^\infty(U)$ we write  $df=\frac{\partial f}{\partial r}dr+d_Hf$. We can express
the vector space $B^1_\Lambda (U)$ of Hamiltonian vector fields as follows:

\begin{eqnarray*}
B^1_\Lambda(U)&=&\{(-1)^{n-1}rX_\Lambda^H(d_Hf_1,\dots,d_Hf_{n-1})\dd r+\\
&+& \sum_{j=1}^{n-1}(-1)^{n-i}r \frac{\partial f_j}{\partial r}X_\Lambda^H(d_Hf_1,\dots,
\widehat{d_Hf_j},\dots.,d_Hf_{n-1})\},
\end{eqnarray*}
with $f_1,\dots,f_{n-1}\in C^\infty(U)$.

Hence all Hamiltonian vector fields must vanish along $H$. For each $H^i$ let us denote  by
$\mathfrak{X}_{free}(H^i)$ the vector space of divergence free vector fields
in $H^i$ with respect to the volume form $\Omega_\Lambda^{H^i}$. Denoting by $r_i$ to the corresponding radial coordinate, we have the following
\begin{corollary}\quad
\begin{enumerate}
\item $\langle r_i\dd {r_i}\rangle \oplus \mathfrak{X}_{free}(H^i)
\subset H^1_\Lambda(U^i)$.
\item ${\bigoplus}_{i \in I}\left(\langle\phi_i\cdot r_i\dd {r_i}\rangle\oplus \phi_i\cdot\mathfrak{X}_{free}(H^i)\right)
\subset  H^1_\Lambda(M)$, where $\phi_i$ are bump functions supported in the collars.
For $n\geq 3$ this space is clearly infinite dimensional.
\end{enumerate}
\end{corollary}

\begin{proof} The assertion about the divergence free vector fields is clear. Regarding the size
of the space we notice that it can be identified with closed $(n-2)$-forms in $H^i$  containing
the exact ones. From the description of $B^1_\Lambda(U)$ we see that the coefficient of $ r\dd {r}$
contains the factor $ X_\Lambda(d_Hf_1,\dots,d_Hf_{n-1})$ which cannot be everywhere non-vanishing 
on each $\{r\}\times H$ by compactness.                           
\end{proof}

We see that the case $n=2$ is quite special and in fact one can easily  compute
$H^1_\Lambda((-1,1)\times S^1)$.

\begin{corollary}  $H^1_\Lambda((-1,1)\times S^1)$ is spanned by the modular vector field
$X_\Lambda^{S^1}$ and $r\dd r $.
\end{corollary}

\begin{proof} The vector field $X_\Lambda^{S^1}$ trivializes $TS^1$ so any vector field can be written 
as $X=A\frac{\partial}{\partial r}+gX_\Lambda^{S^1}$, $A, g\in C^\infty( (-1,1)\times S^1)$. 
One checks that

\begin{equation}\label{eq:cobound} B^1_\Lambda = \{ -rX_\Lambda^{S^1}(d_{S^1}f)\dd r-
r\frac{\partial f}{\partial r}X_\Lambda^{S^1} \;|\; f \in C^\infty((-1,1)\times S^1)\},
\end{equation}
and
\begin{equation}\label{eq:cocycle}
Z^1_\Lambda = \{\left(kr+r\int{\frac{X_\Lambda^{S^1}(d_{S^1}g)}{r}dr} \right)\dd r
+gX_\Lambda^{S^1} \;|\;g \in C^\infty((-1,1)\times S^1),\, g_{\mid {S^1}}=k_1, \,k,k_1 \in \mathbb{R}\}
\end{equation}

And any  $g$ in the above description of a cocycle can be assumed to be the constant $k_1$;
just add the coboundary  defined by  the function $f=\int{\frac{g-k_1}{r}}dr$ 
(as described in \ref{eq:cobound}). 
\end{proof}

Determining the group $H^1_\Lambda(U)$ seems in general a very difficult problem.

Also there seems to be little hope to compute $H^0_\Lambda(M)$ easily. For example for $n=3$ we see
that $X_{f,g}=0$ implies that $df$ and $dg$ have to be proportional. If
we assume $f$ to be a Morse function, $g$ has to be constant on its leaves. So we have
as many choices for $g$ as the ring of smooth functions of the leaf space $M^3/f$. This is a one
dimensional space that  can be very different for the same manifold (one can construct them from a
handle decomposition of the manifold  just looking at how the homotopy group $\pi_0$ changes when we add handles).

\section{Some special families of Nambu structures}

As we have seen, the problem of classifying generic Nambu structures on a given
manifold includes that of the classification of certain arrangements of
connected oriented hypersurfaces (those arrangements that come from the zeros of a function).
For $M^n$ one can consider the dual graph to $(M,\mathcal{H})$ and put a plus sign 
if the orientation of the $n$-tensor in the connected component  coincides with that of $M$,
and minus otherwise. Giving the signs is equivalent to giving the orientation of the hypersurfaces.

For $S^2$, Radko \cite{Ol} defines $\mathcal{G}_k(S^2)$ as the set of generic Poisson structures on $S^2$
with $k$ vanishing curves. She observes that the associated dual graphs are trees (each circle
disconnects the $2$-sphere). A weighted signed tree is defined as a tree with a plus or minus
sign attached to each vertex so that for each vertex, those belonging to the boundary of its star
have opposite sign; each
edge is weighted with a positive number (the modular period), and a real number (the
regularized volume) is assigned to the whole graph. She proves the following:

\begin{theorem}[(\cite{Ol})] The set $\mathcal{G}_k(S^2)$,
up to orientation preserving isomorphisms, coincides with the isomorphism
classes of weighted signed trees  with $k+1$ vertices (the isomorphism has to
preserve the real number attached to the graph).  \end{theorem}

The result  relies on the fact that there is a one to one correspondence between 
arrangements of $k$ circles in $S^2$ (in fact up to isotopy) and isomorphism classes
of trees with $k+1$ vertices (observe also that every tree can be signed in two ways).
One can isotope two arrangements with equivalent tree because, up to isotopy, the
circle sits in $S^2$ in a unique way  splitting $S^2$ in two disks, and that 
results admits a well-known generalization.
 \begin{theorem}[(Smooth Schoenflies theorem)] Any smooth embedding  
$j\colon S^{n-1}\hookrightarrow S^n$ bounds an $n$-dimensional ball and hence splits
the sphere into two $n$-dimensional balls. In particular it is isotopic to the 
standard one where $S^{n-1}$ sits inside
 $\mathbb{R}^n\subset \mathbb{R}^n\cup\{\infty\}=S^n$ as the boundary of the 
Euclidean ball of radius one. It also holds for embeddings in $\mathbb{R}^n$.
 \end{theorem}

 As consequence of this result one easily proves the following:
 \begin{lemma} There is a one to one correspondence between arrangements of $k$
$\;(n-1)$-spheres in $S^n$  and isomorphism classes of trees with $k+1$
vertices. 
 \end{lemma}

\begin{definition} Let us define $\mathcal{G}_k(S^n)$ to be the set of generic 
Nambu structures in $S^n$
whose vanishing set consist of $k$ $(n-1)$-spheres.
\end{definition}
Giving $S^n$ the usual orientation  we can put signs in the dual trees. Thus we have just proved the following
\begin{proposition} The set $\mathcal{G}_k(S^n)$ is, up to isotopy, the same as the equivalence classes of weighted signed trees with $k+1$ vertices (and hence the set of isotopy classes  is the same for every $n\geq 2$).
\end{proposition}


\begin{thebibliography}{xxxxx}

\bibitem{DZ} J.~P.~Dufour and N.~T.~Zung, Linearization of Nambu
    structures, \emph{Compositio Math.~}\textbf{117}, no.~1, 77-98 (1999).

\bibitem{M}  R.~Ib\'a\~nez, M.~de Le\'on, B.~L\'opez, J~C.~Marrero and
 E.~Padr\'on, Duality and modular class of a Nambu structure,
 \emph{J.~Phys.~A} \textbf{34}, no.~17, 3623-3650  (2001).

\bibitem{Mo}  P.~Monnier,  Computations of Nambu-Poisson cohomologies. 
\emph{Int. J. Math. Math. Sci.}  \textbf{26}, no. 2, 65--81.   (2001).
  

\bibitem{Na}  Y.~Nambu, Generalized Hamiltonian mechanics, 
  \emph{Phys.~Rev.~D}, \textbf{7}, 2405-2412 (1973).

\bibitem{Ol} O.~Radko, A classification of topologically stable Poisson structures on a compact oriented surface. \emph{J. Symplectic Geom.}  \textbf{1}, no. 3, 523--542  (2002).

\bibitem{Ta} L.~Takhtajan, On foundations of the generalized Nambu
  mechanics, \emph{Comm.~Math. Phys.~}\textbf{160}, 295-315 (1994).

\bibitem{Va} I.~Vaisman. A survey on Nambu-Poisson brackets,
  \emph{Acta Math.~Univ.~Comenian (N.S.)} \textbf{68}, no. 2, 213-241 (1999).

\end{thebibliography}
\end{document}